\documentclass[12pt,leqno]{article}
\usepackage{amsmath, amsthm, amsfonts, amssymb, color}
\usepackage{mathrsfs}
\newtheorem{thm}{Theorem}[section]

\theoremstyle{definition}

\newcommand{\scr}[1]{\mathscr #1}
\definecolor{wco}{rgb}{0.5,0.2,0.3}

\numberwithin{equation}{section} \theoremstyle{remark}

\newcommand{\ua}{\uparrow}

\title{{\bf Large Deviations for  Stochastic Generalized Porous Media Equations}
\footnote{Supported in part by the DFG through the Forschergruppe
``Spectral Analysis, Asymptotic Distributions and Stochastic
Dynamics'', the BiBoS Research Centre,  NNSFC(10121101) and RFDP.}
}
\author{{\bf  Michael R\"ockner}\\
\footnotesize{Fakult\"at F\"ur Mathematik, Universit\"at
Bielefeld, D-33501 Bielefeld, Germany}\\
{\bf Feng-Yu Wang}\\
\footnotesize{School  of Mathematical Sciences, Beijing Normal
University, Beijing 100875, China}\\
\footnotesize{Email: wangfy@bnu.edu.cn}\\
{\bf Liming Wu}\\
\footnotesize{Laboratoire de Math\'ematiques Appliqu\'ees,}
\\ \footnotesize{CNRS-UMR 6620, Universit\'e Blaise Pascal,
63177 Aubiere, France}}

\begin{document}
\maketitle
\begin{abstract} The large deviation principle is established for the distributions of
 a class of generalized stochastic porous media equations  for both small noise and
short time. \end{abstract} \noindent
 AMS subject Classification:\ 76S05, 60H15.   \\
\noindent
 Keywords: Stochastic porous medium equation, large deviation principle.
 \vskip 2cm

\def\R{\mathbb R}  \def\ff{\frac} \def\ss{\sqrt} \def\BB{\mathbb
B}
\def\N{\mathbb N} \def\kk{\kappa} \def\m{{\bf m}}
\def\dd{\delta} \def\DD{\Delta} \def\vv{\varepsilon} \def\rr{\rho}
\def\<{\langle} \def\>{\rangle} \def\GG{\Gamma} \def\gg{\gamma}
  \def\nn{\nabla} \def\pp{\partial} \def\tt{\tilde}
\def\d{\text{\rm{d}}} \def\bb{\beta} \def\aa{\alpha} \def\D{\scr D}
\def\E{\scr E} \def\si{\sigma} \def\ess{\text{\rm{ess}}}
\def\beg{\begin} \def\beq{\begin{equation}}  \def\F{\scr F}
\def\Ric{\text{\rm{Ric}}} \def\Hess{\text{\rm{Hess}}}\def\B{\mathbb B}
\def\e{\text{\rm{e}}} \def\ua{\underline a} \def\OO{\Omega} \def\b{\mathbf b}
\def\oo{\omega}     \def\tt{\tilde} \def\Ric{\text{\rm{Ric}}}
\def\cut{\text{\rm{cut}}} \def\P{\mathbb P} \def\ifn{I_n(f^{\bigotimes n})}
\def\fff{f(x_1)\dots f(x_n)} \def\ifm{I_m(g^{\bigotimes m})} \def\ee{\varepsilon}
\def\pm{\pi_{{\bf m}}}   \def\p{\mathbf{p}}   \def\ml{\mathbf{L}}
 \def\C{\scr C}      \def\aaa{\mathbf{r}}     \def\r{r}
\def\gap{\text{\rm{gap}}} \def\prr{\pi_{{\bf m},\varrho}}  \def\r{\mathbf r}
\def\Z{\mathbb Z} \def\vrr{\varrho} \def\ll{\lambda}

\section{Introduction and Main Results}

We first recall the existence and uniqueness results on strong
solutions to the stochastic generalized porous media equations
obtained recently in \cite{DRRW}.
 Let $(E,\scr
M,{\bf m})$ be a separable probability space and $(L,\D(L))$ a
negative definite self-adjoint linear operator on $L^2({\bf m})$
with spectrum contained in $(-\infty,-\ll_0]$ for some $\ll_0>0.$

 We
assume that, for a fixed number $r>1$, $L^{-1}$ is bounded in
$L^{r+1}(\m)$, which is e.g. the case if $L$ is a Dirichlet
operator (cf.\ e.g.\ \cite{ma-roe92}) since in this case the
interpolation theorem or simply Jensen's inequality implies
$\|\e^{tL}\|_{r+1}\le \e^{-\ll_0t 2/(r+1)}$ for all $t\ge 0,$
where and in what follows, $\|\cdot\|_p$ denotes the norm in
$L^p(\m)$ for $p\ge 1.$ A classical example of $L$ is the Laplace
operator on a smooth bounded domain in a complete Riemannian
manifold with Dirichlet boundary condition.

Let $H^1:=\D(\ss{-L})$ be the real Hilbert space with  inner
product
$$\<f,g\>_{H^1}:= \<\ss{-L}f, \ss{-L}g\>,$$ where $\<\ ,\ \>$ is
the inner product in $L^2(\m)$. Then the embedding $H^1\subset
L^2(\m)$ is dense and continuous. Let $H:= H^{-1}$ be the dual
Hilbert space of $H^1$ realized through this embedding.

The existence and uniqueness of strong solutions to the following
stochastic differential equation has been proved in  \cite{DRRW}:
\beq\label{**}
  \d X_t = (L\Psi(t,X_t)+\Phi(t,X_t))\d t+Q\d W_t,
\end{equation}
where $Q: L^2(\m)\to H$ is a Hilbert-Schmidt operator with
$q:=q(Q)$ the square of its Hilbert-Schmidt norm, $W_t$ is a
cylindrical Brownian motion on $L^2(\m)$ w.r.t. a complete
filtered probability space $(\OO,\F,\F_t,P)$,

$$\Psi,\Phi: [0,\infty)\times \R\times \OO\to \R$$ are
progressively measurable functions, i.e. for any $t\ge 0$,
restricted on $[0,t]\times \R\times \OO$ they are measurable
w.r.t. $\scr B([0,t])\times\scr B(\R)\times \F_t$, and for any
$(t,\oo)\in [0,\infty)\times \OO,\ \Psi(t,\cdot)(\oo)$ and
$\Phi(t,\cdot)(\oo)$ are continuous on $\R$ and satisfy certain
monotonicity conditions. See \cite{Aronson, AP} for an account of
the classical (deterministic) porous media equations and
\cite{BBDR, BDR, DR1, DR2} for the study of weak solutions and
invariant measures for some stochastic generalized porous media
equations.

To explain what is meant by  strong solutions to (\ref{**}), let
us introduce the embeddings

$$V\subset H\subset V^*$$ as follows.
Consider the reflexive separable Banach space $V := L^{r+1}(\m)$.
Then we can obtain a presentation of its dual space $V^*$ through
the embeddings
$
  V \subset H \equiv H' \subset V^* ,$
where $H$ is identified with its dual through the
Riesz-isomorphism. In other words $V^*$ is just the completion of
$H$ with respect to the norm
\[
  \|f\|_{V^*}
  := \sup_{\|g\|_{r+1}\leq 1} \<f,g\>_H,
  \quad f\in H.
\]
Since $H$ is separable, so is $V^*$. We note that this is
different from the usual representation of $V=L^{r+1}(\m)$ through
the embedding
\[
  V \subset L^2(\m) \equiv L^2(\m)' ,
\]
which, of course, gives $L^{(r+1) / r}(\m)$ as dual. But it is
easy to identify the isomorphism between $L^{(r+1) / r}(\m)$ and
$V^*$. Below we simply use $\<\,,\,\>_H$ to denote
$_{V^*}\<\,,\,\>_{V}$,  i.e. the duality between $V$ and $V^*$,
since $_{V^*}\<\,,\,\>_V = \<\,,\,\>_H$ holds  on $H\times V$. It
is explained in \cite{DRRW} that $L: L^{(r+1)/r}(\m)\to V^*$ is a
densely defined bounded operator, so that it extends uniquely to a
fully defined bounded operator, denoted once again by $L$.
Likewise, the natural embedding $L^2(\m)\subset H\subset V^*$
extends uniquely to a one-to-one map from $L^{(r+1)/r}(\m)$ to
$V^*$ (cf. \cite[Corollary 1.2]{DRRW}). Since $\Psi(t,
v)(\oo),\Phi(t,v)(\oo)\in L^{(r+1)/r}(\m)$, by condition
(\ref{1.1}) below, the map  $b:= L\Psi +\Phi: [0,\infty)\times
V\times \OO\to V^*$ is well-defined.

We assume that there exist two constants $c, \aa>0$ such that

\beq\label{1.1} \beg{split} & |\Psi'(\cdot,s)|+|\Phi'(\cdot,s)| \le c (1+|s|^{r-1}),\\
&_{V}\<u-v, b(\cdot,u)-b(\cdot,v)\>_{V^*}\le
-\aa\|u-v\|_{r+1}^{r+1}+c\|u-v\|_H^2,\ \ \ u,v\in
L^{r+1}(\m)\end{split}\end{equation} holds on $[0,T]\times\OO.$ In
particular, according to \cite{DRRW}, the second inequality in
(\ref{1.1}) holds for some $\aa,c>0$ if there exist constants
$\theta_1> \theta_2/\|L^{-1}\|_{r+1}\ge 0$ and $\si\in \R$ such
that

\beq\label{1.2}\beg{split} &(s-t)(\Psi(\cdot,s)-\Psi(\cdot,t))\ge
\theta_1 |s-t|^{r+1},
\\& |\Phi(\cdot,s)-\Phi(\cdot,t)|\le \theta_2 |s-t|^{r} +\si |s-t|,\ \ \ s,t\in \R\end{split}\end{equation}
holds on $[0,T]\times \OO.$ According to \cite{DRRW} (see also
\cite[Theorems II.2.1 and II.2.2]{KR} for more general
situations), condition (\ref{1.1}) implies that equation
(\ref{**}) has a unique strong solution; that is,
 there is a unique $H$-valued continuous $(\scr F_t)$-adapted process $X_t$
  with $X\in L^{r+1}([0,T]\times \OO\times E, \d t\times P\times
  \m)$ such that for any $e\in L^{r+1}(\m)$,

 \beq\label{D}
    \<X_t,e\>_H
    = \<X_0, e\>_H
      -\int_0^t {\bf m}\big(\Psi(s,X_s)e + \Phi(s,X_s) L^{-1}e\big) \d s
      +\<QW_t,e\>_H,
      \ \ \ t\in [0,T].
  \end{equation}
To see that the solution defined above satisfies the equation
\beq\label{N} X_t = x +\int_0^t(L\Psi +\Phi)(s,X_s)\d s + Q W_t,\
\ \ t\in [0,T]\end{equation} in $H$, we first observe that by
(\ref{1.1}), the right hand side of (\ref{N}) exists in $V^*$ for
any $t>0$ since $X\in L^{r+1}([0,T]\times \OO\times E, \d t\times
P\times \m)$. Since both $X_t-x$ and $QW_t$ take values in $H$,
(\ref{N}) indeed holds in $H$.
\paragraph{Remark 1.1.}  In order to imply the  large deviation principle,
our assumptions are indeed stronger than those used in \cite{KR}
to prove  existence and  uniqueness of strong solutions.  On the
other hand, in \cite{DRRW} we present a direct proof for
existence, uniqueness and ergodicity of strong solutions for
(\ref{**}) under the extra assumption that the spectrum of $L$ is
discrete. Since this assumption was not really used in the proofs,
it can be dropped from that paper. Furthermore, in the recent work
\cite{RRW}, the existence and uniqueness of strong solutions have
been obtained for a much more general framework so that one may
take Orlicz norms in place of $L^{r+1}(\m)$ in applications. Our
arguments  for the large deviation principle presented below are,
however, difficult to be extended to the general situation of
\cite{RRW}.

\ \newline
In this paper we study the large deviation property of
the above stochastic generalized porous medium equation   for both
small noise and short time. Recall (\cite{DZ}) that a sequence of
probability measures $(\mu_\vv)_{\vv>0}$ on some Polish space $E$
satisfies, as $\vv\to 0$, the {\it large deviation principle (LDP
in short) with speed $\lambda(\vv)\to +\infty$ (as $\vv\to 0$) and
rate function} $I: E\to [0,+\infty]$, if $I$ is a good rate
function, i.e., the level sets $\{I\le r\}$, $r\in {\mathbb R}^+$
are compact, and for any Borel subset $A$ of $E$,
$$
 -\inf_{x\in A^o} I(x)\le \liminf_{\vv\to 0}
   \frac 1{\lambda(\vv)} \log \mu_{\vv}(A)\le \limsup_{\vv\to 0}
   \frac 1{\lambda(\vv)} \log \mu_{\vv}(A)\le
  -\inf_{x\in \bar A} I(x),
$$
where $A^o$ and $\bar A$ are respectively the closure and the
interior of $A$ in $E$. In that case we shall simply say  that
$(\mu_\vv)$ satisfies the $LDP(\lambda(\vv), I)$ on $E$, or even
more simply write $(\mu_\vv)\in LDP(\lambda(\vv), I)$ on $E$. We
say that the family of $E$-valued random variables $X^\vv$
satisfies the $LDP(\lambda(\vv), I)$ if the family of their laws
does.

Let us first consider the following stochastic differential
equation with small noise:

\beq\label{1.3} \d X_t^{\vv} = (L\Psi(t,X_t^\vv )+\Phi(t,X_t^\vv
))\d t+\vv Q\d W_t,\ \ \ \vv>0, X_0^\vv =x\in H.\end{equation}

From now on, let  $T>0$ and $x\in H$ be fixed. To state our main
results, let us first introduce the skeleton equation associated
to (\ref{1.3}):

\beq\label{1.4} \ff{\d z_t^\phi }{\d t} = L\Psi(t,z_t^\phi )
+\Phi(t,z_t^\phi ) + \phi_t,\ \ \ z_0^\phi := x,\end{equation}
where $\phi\in L^2([0,T];H).$ An element $z^\phi \in C([0,T]; H)
\cap L^{r+1}([0,T]\times E, \d t \times \m)$ is called a solution
to (\ref{1.4}) if for any $e\in L^{r+1}(\m)$,

\beq\label{1.5} \<z_t^\phi , e\>_H = \<x,e\>_H - \int_0^t\big\{
\<L^{-1} e, \phi_t+\Phi(s,z_s^\phi)\>+ \<e,\Psi(s,z_s^\phi
)\>\big\}\d s,\ \ \ t\in [0,T].\end{equation} We shall prove that
(\ref{1.1}) and (\ref{1.2}) imply the existence and the uniqueness
of the solution to (\ref{1.4}) for any $\phi\in L^2([0,T];H)$, and
thus, as explained above for the solution to $(\ref{**})$, the
solution satisfies the corresponding integral equation of
(\ref{1.4}) in $H$.

Now, we introduce the rate function. For any $\phi\in
L^2([0,T]\times E,\d t\times \m),$ let $\|\phi\|_{L^2}^2:=
\int_0^T\d t\int_E \phi_t^2\d\m.$ Define

\beq\label{rate} I (z):= \ff 1 2 \inf \{\|\phi\|_{L^2}^2:\
z=z^{Q\phi} ,\ \phi\in L^2([0,T]\times E,\d t\times \m)\},\ \ \
z\in C([0,T];H), \end{equation} where we set $\inf\emptyset
=\infty$ by convention. The following result is of a
Freidlin-Wentzell type estimate:

\beg{thm}\label{T1.1} Assume $(\ref{1.1})$. For each $\vv>0$, let
$X^\vv=(X^\vv_t)_{t\in [0,T]}$ be the solution to $(\ref{1.3})$.
Then as $\vv\to 0$, $(X^\vv)$ satisfies the $LDP(\vv^{-2}, I)$ on
$C([0,T]; H) ($equipped with the sup-norm topology$)$, where the
rate function $I$ is given by $(\ref{rate})$.
\end{thm}

Next, we consider the LDP of the solution $X_t$ to $(\ref{**})$
for short time, which in the classical finite dimensional case is
the famous Varadhan's large deviation estimate. Since $X_{\vv^2
t}$ solves the
 equation

\beq\label{st} \d \tt X_t^{\vv} = \vv^2(L\Psi(\vv^2t,\tt X_t^\vv
)+\Phi(\vv^2t,\tt X_t^\vv ))\d t+\vv Q\d \tilde W_t,\ \ \ \vv>0,
\tt X_0^\vv =x,\end{equation} where $(\tilde
W_t:=(1/\vv)W_{\vv^2t})$ is a BM of the same law as $(W_t)$, it
suffices to establish the LDP for the law of $\tt X^\vv$.

\beg{thm}\label{T1.2} Assume $(\ref{1.1})$. If $x\in L^{r+1}(\m)$
then $\tt X^\vv=(X_{\vv^2 t})$ satisfies the $LDP(\vv^{-2}, \tt
I)$ where
$$
\tt I(z):=  \ff 1 2\inf\bigg\{  \|\phi\|_{L^2}^2:\ z_t=
x+Q\int_0^t \phi_s\d s\bigg\},\ \ z\in C([0,T];H).
$$
\end{thm}

Let us make some historical comments. In the finite dimensional
case, under the Lipschitzian condition, the LDP of
$X_{\vv^2\cdot}$ is the famous Varadhan's estimate \cite{Var}, and
Theorem \ref{T1.1} is the well known Freidlin-Wentzell's LDP
(\cite{FW}). For the extensions to infinite dimensional diffusions
or stochastic PDE under global Lipschitz condition on the
nonlinear term, we refer the reader to Da Prato and Zabczyk
\cite{DaZa}(also for the literature until 1992). For the case of
local Lipschitz conditions we refer to \cite{CR} where also
multiplicative and degenerate noise is handled. Unlike in our
situation, in \cite{CR} the drift still contains a nontrivial
(therefore smoothing) linear part. In many examples of SPDE,
however, (local) Lipschitz conditions are rarely satisfied (such
as the porous equation in this work). Without Lipschitz
conditions, each type of stochastic non-linear PDE requires
specific techniques and adapted estimates. So the situation
becomes much more dispersive. Here we mention only  the work of
Cardon-Weber \cite{C-W} on the LDP for stochastic Burgers
equations with small noise and the important work of Hino and
Ramirez \cite{Hino} for the Varadhan's small time estimate of
large deviations for general symmetric Markov processes, where the
reader may also find other recent references.

Here are some remarks on Theorem \ref{T1.2} related with  the
general work of Hino and Ramirez \cite{Hino}: 1) As our process
$(X_t)$ is highly non-symmetric,  the result in \cite{Hino} can
not be applied. 2) The extra condition on $x\in L^{r+1}(\m)$ (not
all $x\in H$) in Theorem \ref{T1.2} is also a quite general
phenomenon in infinite dimension because the result of \cite{Hino}
holds only for $\mu-a.e. x$ where $\mu$ is the invariant measure,
and in our case, the invariant measure is supported in
$L^{r+1}(\m)$ (\cite{DRRW}). 3) Furthermore the LDP in Theorem
\ref{T1.2} is pathwise, unlike that in \cite{Hino} which is only
for the marginal law.

This paper is organized as follows. The next section is devoted to
the study of the skeleton process $z^\phi$, which is crucial for
identifying the rate function of our LDP. In \S 3 we give an {\it
a priori} exponential estimate and recall the generalized
contraction principle. The proof of Theorem \ref{T1.1} is
presented in \S 4, and our strategy is based on two procedures of
approximation: first for finite dimensional noise (i.e., only a
finite number of directions are stochastically perturbed) we
approximate the path of $QW$ piecewise linear; second, we
approximate the whole noise $QW$ by the finite dimensional noises.
This strategy can be easily adapted for the proof of Theorem
\ref{T1.2} in \S 5.

\section{The skeleton process}

\beg{prp} \label{P2.2} Assume $(\ref{1.1})$. Let
$\|z\|:=\sup_{t\in [0,T]}\|z_t\|_H$ for $z\in C([0,T];H).$
 For any $x\in H$ and any $\phi\in L^2([0,T];H)$ there
 exists a unique solution $z^\phi $ to $(\ref{1.4})$ and

 \beq\label{NN}\int_0^T \m(|z_t^\phi-z_t^\psi |^{r+1})\d t \le
 C\int_0^T\|\phi_t-\psi_t\|_H^2\d t,\end{equation}
 \beq\label{stable} \|z^\phi -z^{\psi} \|\le C\int_0^T
\|\phi_t-\psi_t\|_H\d t\end{equation}
 hold for some constant $C>0$ and all $x\in H,\ \phi,\psi\in
 L^2([0,T];H).$
\end{prp}

\beg{proof} To verify the existence of the solution, we  make use
of \cite[Theorem II.2.1]{KR}.
 Let $V:= L^{r+1}(\m)$
and $V^*$ the duality of $V$ w.r.t. $H$, and let $B:=0$ and

$$A(s,v):= L\Psi(s,v)+\Phi(s,v) + \phi_s.$$
Then, due to (\ref{1.1}),  it is trivial to verify Assumptions
$A_i) (i=1,..,5)$ on page 1252 of \cite{KR} for some $K,\aa>0,
p:=r+1, q:=\ff{r+1} r, $ and $f(t):= c(1+\|\phi_t\|_H^q)$ for some
constant $c>0.$ Then, by \cite[Theorems II.2.1 and II.2.2]{KR}
(see also \cite[Theorem 30.A]{Z}) (\ref{1.4}) has a unique
solution. Let $z^\phi$ be the unique solution to (\ref{1.4}) for
$\phi\in L^2([0,T];H).$

By It\^o's formula due to \cite[Theorem I.3.2]{KR} and
(\ref{1.1}), we have

\beq\label{RW}\beg{split}\ff{\d}{\d t} \|z_t^\phi -z_t^\psi\|_H^2
= &-2 \<z_t^\phi-z_t^\psi, \Psi(t,z_t^\phi )
-\Psi(t,z_t^\psi)\>\\
&- 2 \<L^{-1}(z_t^\phi-z_t^\psi) ,
\Phi(t,z_t^\phi )-\Phi(t,z_t^\psi) +\phi_t-\psi_t \>\\
\le &-2\aa \m(|z_t^\phi-z_t^\psi |^{r+1})+2c\|z_t^\phi
-z_t^\psi\|_H^2
+2\|z_t^\phi-z_t^\psi\|_H\|\phi_t-\psi_t\|_H.\end{split}\end{equation}
 Since (\ref{RW})  implies, for any
$\vv>0$, that

$$\ff{\d}{\d t} (\vv+\|z_t^\phi -z_t^\psi\|_H^2)^{1/2}\le
\|\phi_t-\psi_t\|_H +c \|z_t^\phi-z_t^\psi\|_H,$$ by Gronwall's
lemma we have

$$\e^{-cT}\ss{\vv+\|z^\phi-z^\psi\|^2} \le \vv
+\int_0^T\|\phi_t-\psi_t\|_H\d t.$$ This implies (\ref{stable})
for $C:=\e^{c T}$ by letting $\vv\to 0.$  Finally, (\ref{NN})
follows by combining (\ref{stable}) with (\ref{RW}).
\end{proof}

\section{Exponential estimates and a generalized contraction principle}

The following {\it a priori} estimate will be crucial for the
proof of Theorem \ref{T1.1}.

\beg{lem}\label{L3.1} Assume $(\ref{1.1})$. Then for any
$\gg>0,q_0>0$ and $\vv_0>0$ there exits a constant $c>0$ such that
for all $Q$ with $q(Q)\le q_0$ and all $\vv\in (0,\vv_0)$,

\beq\label{3-1} \mathbb E \exp\left(\gg \vv^{-2}\int_0^T
\|X^\vv_t\|_{r+1}^{r+1}\d t\right)\le e^{c\vv^{-2}}.
\end{equation}
\end{lem}
Throughout this paper we adopt the following notation: for two
continuous real semimartingales $(x_t)$ and $(y_t)$, $dx_t\le
dy_t$ means that their martingale parts are the same and
$x_t-x_s\le y_t-y_s$ for all $t>s\ge 0$.

\beg{proof} By $(\ref{1.1})$, $(\ref{1.2})$ with $\theta_2<
\theta_1$ and using It\^o's formula due to \cite[Theorem
I.3.2]{KR}, there exist constants $c_0, c_1, c_2>0$ and $\vv_0>0$
such that for all $\vv<\vv_0$,
$$\aligned
d\|X_t^\vv\|_H^2 &\le  -2\<X_t^\vv, (\Psi
+L^{-1}\Phi)(t,X_t^\vv)\> \d t+  2\vv
\<X_t^{\vv}, Q\d W_t\> + q\vv^2\d t\\
&\le - {\bf m}\left(\aa
|X_t^\vv|^{r+1}-c_0[|X_t^\vv|^{2}+1]\right)
\d t+ 2\vv \<X_t^{\vv}, Q\d W_t\> +(q \vv^2+c_0\|X_t^\vv\|_H^2)\d t\\
&\le -c_1{\bf m}(|X_t^{\vv}|^{r+1})\d t+ c_2 \d t+2\vv
\<X_t^{\vv}, Q\d W_t\>_H.
\endaligned
$$
Then \beq\label{b31} \|X_T^\vv\|_H^2 -\|x\|_H^2 + c_1 \int_0^T{\bf
m}(|X_t^{\vv}|^{r+1}) \d t \le 2\vv \int_0^T\<X_t^{\vv}, Q\d
W_t\>_H +c_2T.
\end{equation}
Letting $dM_t:= \<X_t^{\vv}, Q\d W_t\>_H$ (with $M_0=0$), since
$\forall \lambda\in \R, \ \xi_t:=\exp(\lambda M_t-\frac
{\lambda^2}2\<M\>_t)$ is a martingale and the quadratic
variational process $\<M\>_t$ satisfies $d\<M\>_t\le
q_0\|X_t^\vv\|_H^2\d t$, we obtain from (\ref{b31}) that, for
$\ll:= 8\gg/c_1\vv,$

\beg{equation*}\beg{split} &\mathbb E
\exp\bigg(\gg\vv^{-2}\int_0^T\|X_t^\vv\|_{r+1}^{r+1}\d
t\bigg)=\mathbb E
\exp\bigg(2\gg\vv^{-2}\int_0^T\|X_t^\vv\|_{r+1}^{r+1}\d t
-\gg\vv^{-2}\int_0^T\|X_t^\vv\|_{r+1}^{r+1}\d t\bigg)\\
&\le \mathbb E\exp\bigg(\ff {4\gg}{c_1\vv}M_T +
\ff{2c_2\gg}{c_1\vv^2}T +\ff{2\gg\|x\|_H^2}{c_1\vv^2} -\gg \vv^2
\int_0^T\|X_t^\vv\|_{r+1}^{r+1}\d t\bigg)\\
&\le \mathbb E\exp\bigg(\ff \ll 2 M_T -\ff{\ll^2} 4 \<M\>_T
+\int_0^T\Big(
\ff{q_0\ll^2}4\|X_t^\vv\|_H^2-\gg\vv^{-2}\|X_t^\vv\|_{r+1}^{r+1}\Big)\d
t +\ff{2c_2\gg}{c_1\vv^2}T +\ff{2\gg\|x\|_H^2}{c_1\vv^2}\bigg)\\
&\le \Big\{\mathbb E \xi_T\Big\}^{1/2}\bigg\{\mathbb
E\exp\bigg(\int_0^T\Big(
\ff{q_0\ll^2}2\|X_t^\vv\|_H^2-2\gg\vv^{-2}\|X_t^\vv\|_{r+1}^{r+1}\Big)\d
t +\ff{4c_2\gg}{c_1\vv^2}T
+\ff{4\gg\|x\|_H^2}{c_1\vv^2}\bigg)\bigg\}^{1/2}\\
&\le \exp(c\vv^{-2}T)\end{split}\end{equation*} for some constant
$c>0$ and all $\vv\in (0,\vv_0)$, where the last step is due to
the martingale property of $\xi_t$ and that $\|\cdot\|_{r+1}\ge
c\|\cdot\|_H$ for some $c>0$ and that $r>1.$\end{proof}

In large deviation theory, when $(\mu_\vv)$ satisfies the
$LDP(\lambda(\vv), I)$ on a Polish space ${\bf E}$ and if $f: {\bf
E}\to {\bf F}$ is continuous where ${\bf F}$ is another Polish
space, then $(\mu_\vv\circ f^{-1})\in LDP(\lambda(\vv), I_f)$,
where
$$
I_f(z):=\inf_{f^{-1}(z)}I, \ \ \ \  z\in {\bf F}.
$$
That is the so called contraction principle. The following
generalization is taken from \cite{Wu04} (some preceding weaker
versions can be found in \cite[Theorems 4.2.16 and 4.2.23]{DZ}).

\begin{thm}\label{thma31}
{\bf (Generalized Contraction Principle)} Let ${\bf E}, {\bf F}$
be two Polish spaces and $(\mu_\vv)$ a family of probability
measures on $E$. If $(\mu_\vv)\in LDP(\lambda(\vv), I)$ and there
exists a sequence of continuous mappings $f^N : {\bf E}\to {\bf
F}$ such that

\begin{equation}\label{145}
\lim_{N\to\infty}\limsup_{\vv\to0} \frac 1{\lambda(\vv)} \log
\mu_\vv \left(\rho_F\left(f^N, f\right)>\delta \right)=-\infty, \
\ \ \delta >0,
\end{equation}
where $\rho_F$ is some compatible metric on ${\bf F}$ and $f: {\bf
E}\to {\bf F}$ is a measurable mapping, then there exists a
continuous function $\tilde f : \{I<+\infty\}\to {\bf F}$ such
that

\begin{equation}
\lim_{N\to\infty}\sup_{I\le r}\rho_F(f^N, \tilde f) =0, \ \ \ \
r>0; \label{146}
\end{equation}
and $(\mu_\vv(f\in\cdot))\in LDP(\lambda(\vv), I_{\tilde f})$,
where

\begin{equation}
I_{\tilde f}(z):=\inf_{\tilde f^{-1}(z)} I,\ \ \ \  z\in {\bf F}.
\label{147}
\end{equation}
\end{thm}

\section{Proof of  Theorem \ref{T1.1}}

We shall prove Theorem \ref{T1.1} by two procedures of
approximation. Let $\{e_i:i\ge 1\}$ be dense in $L^{r+1}$ and
hence, also dense in $H$.
 For any fixed $n\ge 1$, let
$H_n:=\text{span}\{e_i: 1\le i\le n\}$ and $P_n:H\to H_n$ be the
orthogonal projection. Let  $X_{t}^{\vv,n}$ be the solution of

\beq\label{finite} \d X_{t}^{\vv,n}= (L\Psi+\Phi)(t,X_t^{\vv,n})\d
t +\vv P_n Q \d W_t,\ \ \ X_0^{\vv, n}=x.\end{equation}

Next for each $N\in\N$ and for any path $w\in C([0,T];H)$, let
$t_i:= iT/N$ for $0\le i\le N$ and define the ($N-$times)
piecewise linear approximation of $w$ by

$$w_{t}^{(N)}:= \ff N T\sum_{i=0}^{N-1} 1_{(t_i,t_{i+1}]}(t) \big(
(t-t_{i}) w_{t_{i+1}} + (t_{i+1}-t) w_{t_i}\big),\ \ \ t\in [0,T].
$$
By Proposition \ref{P2.2}, the following equation has a unique
solution $X_{t,N}^{\vv,n}$ in $H$:

\beq\label{tN} \dot X_{t,N}^{\vv,n}:=\ff{\d X_{t,N}^{\vv,n}}{\d
t}= L\Psi(t,X_{t,N}^{\vv,n}) +\Phi(t,X_{t,N}^{\vv,n})+ \vv \frac
d{\d t} (P_nQ W)_{t}^{(N)},\ \ \ X_{0,N}^{\vv,n}= x.\end{equation}

We claim that it is enough to establish

\beq\label{prob1} \limsup_{N\to\infty}\limsup_{\vv\to 0}\vv^2\log
\mathbb P (\|X^{\vv,n}- X_{\cdot,N}^{\vv,n}\|>\dd)=-\infty,\ \
\forall \dd>0\end{equation} and

\beq\label{W2} \limsup_{n\to\infty}\limsup_{\vv\to 0}\vv^2\log
\mathbb P (\|X^{\vv,n}- X^{\vv}\|>\dd)=-\infty,\ \ \forall
\dd>0.\end{equation}

In fact, by  Schilder's theorem, the law of $\vv QW$ satisfies the
LDP on $C([0,T]; H)$ with speed $\lambda(\vv)=\vv^{-2}$ and with
rate function given by
$$
J(\tilde \phi)=
   \inf\left\{\frac 12 \|\phi\|_{L^2}^2: \
    \tilde Q\phi=\tilde \phi(\cdot)\right\}
$$
where $\tilde Q : \phi\mapsto \int_0^{\cdot} Q \phi(s) \d s$ is a
continuous linear mapping from  $L^2([0,T]\times E, \d t\times
\m)$ to $C([0,T];H)$ (with the convention that $\inf\emptyset
:=+\infty$). Next, let $f_{n,N}$ denote the map which associates
each path $\omega\in C([0,T]; H)$ of $\vv QW$ to the solution
$X_{t,N}^{\vv,n}$ of (\ref{tN}), i.e., $\gamma:=f_{n,N}(\omega)$
is the the unique solution of
$$
\gamma_t = x+ \int_0^t [L\Psi(s,\gamma_s)+\Phi(s,\gamma_s)] \d s
+(P_nw)_{t}^{(N)},
$$
where $(P_nw)_{t}^{(N)}$ is the ($N-$times) piecewise linear
approximation of $P_nw$. Applying Proposition \ref{P2.2} with $Q$
replaced by $P_nQ$ and noting that
$(P_nw)_{t}^{(N)}=P_n(P_nw)_{t}^{(N)}$ and that all norms on $H_n$
are equivalent, we see that $f_{n,N}: C([0,T]; H) \to C([0,T]; H)$
is continuous. Furthermore, by (\ref{prob1}) and (\ref{W2}), for
each $n$, there is some $N(n)$ such that
$$
\limsup_{n\to\infty}\limsup_{\vv\to 0}\vv^2\log \mathbb P
(\|X^{\vv}- f_n(\vv QW)\|>\dd)=-\infty,\ \ \forall \dd>0
$$
where $f_n:=f_{n,N(n)}$. Hence by Theorem \ref{thma31}, $X^{\vv}$
satisfies the LDP on $C([0,T]; H)$ with rate function given by
$$
I(z)=\inf\{J(\tilde \phi):\ \tilde f(\tilde \phi)=z\} =\inf\{\frac
12 \|\phi\|_{L^2}^2:\ \tilde f(\tilde Q\phi) =z\}
$$
and $\tilde f(\tilde Q\phi)=\lim_{n\to\infty} f_{n, N(n)}(\tilde
Q\phi)$ by (\ref{146}). But by Proposition \ref{P2.2} and the
following Lemma \ref{L4.1}, $f_{n, N}(\tilde Q\phi)\to z^{\phi}$
as $n,N$ goes to infinity. Thus $\tilde f(\tilde Q\phi)=z^{\phi}$,
which yields the claimed rate function.

\beg{lem}\label{L4.1} For any $\phi\in L^2([0,T]\times E,\d
t\times\m)$, let $h_t:= \int_0^t Q\phi_s\d s, t\in [0,T].$ For any
sequence $N(n)\to\infty$ as $n\to\infty$ we have

$$\lim_{n\to\infty} \int_0^T \Big\|\ff{\d}{\d t} (P_n h)_t^{N(n)}
-Q\phi_t\Big\|_H^2\d t=0.$$\end{lem}

\beg{proof} Let $t_i:= Ti/N(n)$. Since

$$\ff{\d}{\d t} h^{(N(n))} = \sum_{i=1}^{N(n)}
1_{[t_{i-1},t_i)}\ff {N(n)} T\int_{t_{i-1}}^{t_i} Q\phi_s \d s$$
which converges to $Q\phi$ in $L^2([0,T];H)$ as $n\to\infty$, it
suffices to prove that

\beq\label{L2}I_n := \int_0^T \Big\|\ff{\d}{\d t} (P_n
h)_t^{(N(n)} -\ff{\d}{\d t} h_t^{(N(n))}\Big\|_H^2\d t\to
0\end{equation} as $n\to\infty.$ Note that for any $\psi\in
L^2([0,T];H)$ we have

\beg{equation*}\beg{split}&\int_0^T \Big\|\ff{\d}{\d t}
\Big(\int_0^t \psi_s\d s\Big)^{(N(n))}\Big\|_H^2\d t=
\sum_{i=1}^{N(n)}\ff T{N(n)} \Big\|\ff {N(n)}
T\int_{t_{i-1}}^{t_i} \psi_t\d t\Big\|_H^2\\ &\le
\sum_{i=1}^{N(n)}\int_{t_{i-1}}^{t_i} \|\psi_t\|_H^2\d
t=\int_0^T\|\psi_t\|_H^2\d t.\end{split}\end{equation*} Then

$$\lim_{n\to\infty}I_n\le \lim_{n\to\infty}\int_0^T \|P_n Q\phi_t- Q\phi_t\|_H^2\d t=0.$$
\end{proof}

So, to finish the proof of Theorem \ref{T1.1}, we have to prove
(\ref{prob1}) and (\ref{W2}) which will be done in the following
two subsections.

\subsection{ Proof of (\ref{prob1})}

Let $b:= L\Psi+ \Phi$, and

$$
\aligned \hat X_t&:= X_t^{\vv,n} -X_{t,N}^{\vv,n},\ \ \hat\gg_t:=
\vv P_nQ (W_t - W_{t}^{(N)}),
\endaligned
$$
By
(\ref{finite}) and (\ref{tN}) we have

$$\ff{\d \|\hat X_t -\hat \gg_t\|_H^2}{\d t} = 2\<\hat X_t -\hat
\gg_t, b(t,X_t^{\vv,n} -\hat \gg_t)- b(t,X_{t,N}^{\vv,n})\>_H +
2\<\hat X_t -\hat \gg_t, b(t,X_{t}^{\vv,n})
-b(t,X_t^{\vv,n}-\hat\gg_t)\>_H.$$ Combining this with (\ref{1.1})
and (\ref{1.2}) with $\theta_2 <\theta_1$, and using Young's
inequality $xy\le x^{r+1}/(r+1) + [r/(r+1)]y^{(r+1)/r}, \forall
x,y\ge0$, we conclude that there exist $\ll,c>0$ and
$c(\lambda)>0$ such that

\beq\label{1} \aligned &\ff{\d \|\hat X_t -\hat
\gg_t\|_H^2\e^{-ct}}{\d t}\\ &\le -\ll \|\hat X_t -\hat
\gg_t\|_{r+1}^{r+1}  + 2 \|\hat X_t -\hat \gg_t\|_{r+1}\| L^{-1}
(b(t,X_t^{\vv,n})-
b(t,X_t^{\vv,n}-\hat\gg_t))\|_{(r+1)/r}\\
&\le c(\lambda)\| L^{-1} (b(t,X_t^{\vv,n})-
b(t,X_t^{\vv,n}-\hat\gg_t))\|_{(r+1)/r}^{(r+1)/r}+c\|\hat X_t
-\hat \gg_t\|_H^2. \endaligned
\end{equation}
Since $|\Psi'(s)|+|\Phi'(s)|\le c(1+ |s|^{r-1})$ and $L^{-1}$ is
bounded in $L^{(r+1)/r}(\m)$, there exist $c_1, c_2>0$ such that

\beg{equation}\label{2} \aligned &\| L^{-1} (b(t,X_t^{\vv,n})-
b(t,X_t^{\vv,n}-\hat\gg_t))\|_{(r+1)/r}^{(r+1)/r} \le c_1 \|
|\hat\gg_t|(1 + |\hat \gg_t|^{r-1}
+|X_t^{\vv,n}|^{r-1})\|_{(r+1)/r}^{(r+1)/r}\\
&\le c_2 \int_E(|\hat \gg_t|^{r+1} +|\hat \gg_t|^{(r+1)/r} + |\hat
\gg_t|^{(r+1)/r} |X_t^{\vv,n}|^{(r^2-1)/r}) \d{\bf m}.
\endaligned
\end{equation}
 From (\ref{1}) and (\ref{2}) and Young's inequality we obtain that for each $R>1$,

\beg{equation}\label{a31}\beg{split}\ff{\d \|\hat X_t -\hat
\gg_t\|_H^2\e^{-ct}}{\d t} \le c_2 \big\{(1+R)\| \hat
\gg_t\|_{r+1}^{r+1} + \|\hat \gg_t\|_{(r+1)/r}^{(r+1)/r} +
c(r)R^{-1/(r-1)}\|X_t^{\vv,n}\|_{r+1}^{r+1}\big\}\
\end{split}\end{equation}
for some $c(r)>0$. Since all $L^p$-norms ($1\le p\le r+1$) on
$H_n$ are equivalent, for any norm $\|\cdot\|_p$ on $H_n$,  by the
LDP of $\vv P_nQ W_t$ on $C([0,T]; H_n),$ whose good rate function
is denoted by $I_n,$ and recalling that $\hat\gg_t:= \vv
P_n[(QW)_t - (QW)_{t}^{(N)}]$, we have

$$
\aligned
&\limsup_{N\to\infty} \limsup_{\vv\to 0}\vv^{-2}
\log\mathbb P \Big( \sup_{t\in [0,T]}\|\hat \gg_t\|_p>\dd\Big)\\
&\le \limsup_{N\to\infty} -\inf\{I_n(w):\ w\in C([0,T]; H_n),
\sup_{t\in [0,T]}\|w_t - w_{t,N}\|_p\ge \dd \}=-\infty,\ \ \
\forall\dd>0,
\endaligned
$$
where the equality follows from the fact that $\inf_{F_N}I_n\to
+\infty\ (N\to\infty)$ for any sequence of closed subsets
decreasing to $\emptyset$ (an elementary property of a good rate
function).

 Combining this with (\ref{a31}), we see that for any $\dd\in (0,1)$, there exists $c_3>0$ such that the
 l.h.s. of (\ref{prob1}) is less  than
 $$
\limsup_{n\to\infty}\limsup_{\vv\to 0}\vv^2\log \mathbb P
\left(c_3R^{-1/(r-1)} \int_0^T \|X_t^{\vv,n}\|_{r+1}^{r+1} \d t
>\dd\right)
 $$
 which goes to $-\infty$ when $R\to +\infty$ by Chebychev's inequality and the a priori
 exponential estimate in Lemma \ref{L3.1}.

 \subsection{Proof of (\ref{W2})}

  By
(\ref{1.1})and using It\^o's formula in \cite[Theorem I.3.2]{KR},
we have

$$\d \|X_t^\vv -X_t^{\vv,n}\|_H^2 \le (\dd(n)
\vv^2+c\|X_t^\vv -X_t^{\vv,n}\|_H^2)\d t + 2\vv \d M_t^{(n)},$$
where $c>0$ is a constant, $\dd(n):=q(P_nQ-Q)$ is the square of
the Hilbert-Schmidt norm of $P_n Q- Q$ from $L^2(\m)$ to $H$, and
$\d M_t^{(n)}:= \<X_t^\vv-X_t^{\vv,n}, (I-P_n)Q\d W_t)\>_H.$ The
quadratic variation process of the local martingale $M^{(n)}$
verifies
$$
d\<M^{(n)}\>_t\le \|X_t^\vv-X_t^{\vv,n}\|_H^2 \delta(n) \d t.
$$
For any constant $\aa>0$, let $\xi_t:= \exp[\aa \vv^{-2}
\|X_t^\vv-X_t^{\vv,n}\|_H^2\e^{-(1+c)t}]=:\exp[\aa \vv^{-2}Y_t]$.
We have by It\^o's formula in \cite[Theorem I.3.2]{KR} that

$$
\aligned \d \xi_t &\le 2\aa\vv^{-1}\e^{-(1+c)t} \xi_t\d
M_t^{(n)}\\
&\qquad + \vv^{-2} \aa\e^{-(1+c)t}\xi_t\big\{ \dd(n)\vv^2 -
\|X_t^\vv-X_t^{\vv,n}\|_H^2+ 2\aa\e^{-t} \dd(n)\|X_t^\vv-
X_t^{\vv,n}\|_H^2\big\}\d t\\
&\le 2\aa\vv^{-1}\e^{-(1+c)t} \xi_t\d M_t^{(n)} + \aa\dd(n) \xi_t
\d t,
\endaligned
$$
once $1\ge 2\aa\dd(n)$ which holds for all sufficiently large $n$
for $\dd(n)\to 0$ as $n\to\infty$. So $N_t:= \xi_t \exp[-\aa
\dd(n)t]$ is a supermartingale. Therefore, for all $n$ large
enough,

\beg{equation*}\beg{split}\mathbb P(\|X^\vv-X^{\vv,n}\|>\dd)&\le
\mathbb P\Big(\sup_{t\in [0,T]}N_t
>\exp[\dd^2\aa\vv^{-2}\e^{-(1+c)T}-\aa\dd(n)T]\Big)\\
&\le \exp[-\aa\dd^2 \vv^{-2}\e^{-(1+c)T}+\aa\dd(n)
T].\end{split}\end{equation*} This implies (\ref{W2}) since
$\aa>0$ was arbitrary.

\section{Proof  of Theorem \ref{T1.2}}

\ \newline\emph{Proof of Theorem \ref{T1.2}.} (a) We first assume
that there exists   $n\in\N$ such that $q_{ij}=0$ for $i>n$. In
this case the law of $\vv Q W_t +x$ satisfies the large deviation
principle with the given rate function of compact level sets. Thus
by the approximation lemma in large deviations (see \cite[Theorem
4.2.13]{DZ}), it suffices to show that

\beq\label{W3}\limsup_{\vv\to 0}\vv^2\log \mathbb P (\|\tt
X^{\vv}-x-\vv Q\tt W\|>\dd)=-\infty,\ \ \ \dd>0.\end{equation} By
(\ref{1.1}) and (\ref{1.2}) with $\theta_2< \theta_1$, there
exists $\ll,c,c_0>0$ such that

\beg{equation*}\beg{split}\ff{\d \|\tt X_t^\vv -\vv Q\tt
W_t-x\|_H^2\e^{-ct}}{\d
t} &\le -\ll\vv^2 \|\tt X_t^\vv -\vv Q\tt W_t-x\|_{r+1}^{r+1} \\
&\ \ \ + 2\vv^2\e^{-ct}\<
\tt X_t^\vv -\vv Q\tt W_t-x, (L\Psi +\Phi)(\vv Q\tt W_t +x)\>_H\\
&\le c_0\vv^2 \|(\Psi+L^{-1}\Phi)(x+\vv Q\tt
W_t)\|_{(r+1)/r}^{(r+1)/r},\ \ \ t\in
[0,T].\end{split}\end{equation*} Since $L^{-1}$ is bounded in
$L^{(r+1)/r}(E,\m)$, $|\Psi'(s)|+|\Phi'(s)|\le c(1+|s|^{r-1})$ for
some $c>0$, and $x\in L^{r+1}$, there exists $c_1>0$ such that

$$\ff{\d \|\tt X_t^\vv -\vv Q\tt W_t-x\|_H^2\e^{-ct}}{\d t}\le
c_1\vv^2(\|\vv Q\tt W_t\|_{r+1}^{r+1}+1). $$ This immediately
implies (\ref{W3}) by the LDP of $\vv Q\tt W$ in $C([0,T]; H_n)$.
Note that on $H_n$ the norms $\|\cdot\|_H$ and $\|\cdot\|_{r+1}$
are equivalent.

(b) In general, for any $n\ge 1$ let $Q^{(n)}:= P_nQ.$ By (a), the
law of $\tt X_t^{\vv,n}$, the solution to (\ref{st}) for $Q^{(n)}$
in place of $Q$, satisfies the LDP with the good rate function

$$\tt I_n(z) := \ff 1 2\inf\bigg\{  \|\phi\|_2^2:\ z_t= x+\int_0^t
Q^{(n)}\phi_s\d s\bigg\},\ \ z\in C([0,T];H).$$ Similarly to  the
proof of (\ref{prob1}) in \S 4.2  we have

$$\limsup_{n\to\infty}\limsup_{\vv\to 0}\vv^2\log \mathbb P
(\|\tt X^{\vv,n}-\tt  X^{\vv}\|>\dd)=-\infty,\ \ \ \dd>0.$$
Moreover, since $\dd_n:= \|L^{-1/2}(Q-Q^{(n)})\|_{2\to 2} \to 0$
as $n\to\infty$ and since

$$\int_0^T \|(Q-Q^{(n)})\phi_t\|_H\d t\le
\dd(n)\int_0^T\|\phi_t\|_2\d t\le \dd(n) \ss T \|\phi\|_2,$$ we
conclude that the law of $\tt X^\vv$ satisfies the LDP with the
claimed rate function $I$ by the approximation lemma (see
\cite[Theorem 4.2.13]{DZ}). \qed

\beg{thebibliography}{99}

\bibitem{Aronson} D.G. Aronson, \emph{The porous medium equation,} Lecture Notes Math.
Vol. 1224, Springer, Berlin, 1--46, 1986.

\bibitem{AP} D.G. Aronson and L.A. Peletier, \emph{Large time behaviour of solutions
of the porous medium equation in bounded domains,} J. Diff. Equ.
39(1981), 378--412.

\bibitem{BBDR} V. Barbu, V.I. Bogachev,  G. Da Prato and M. R\"ockner, \emph{Weak
solution to the stochastic porous medium  equations: the
degenerate case,} preprint.

\bibitem{BDR} V.I. Bogachev,  G. Da Prato and M. R\"ockner, \emph{Invariant measures
of   stochastic generalized porous medium equations,} to appear in
Dokl. Math.

\bibitem{C-W} C. Cardon-Weber, \emph{Large deviations for a Burger's type SPDE,}
 Stoch. Proc. Appl. 84(1999), 53-70.

\bibitem{CR} S. Cerrai and M. R\"ockner, \emph{Large deviations
for stochastic reaction-diffusion systems with multiplicative
noise and non-Lipschitz reaction term,} Ann.  Probab. 32(2004),
1100--1139.

\bibitem{DR1} G. Da Prato and M. R\"ockner, \emph{Weak solutions to stochastic
porous media equations,} J. Evolution Equ. 4(2004), 249--271.

\bibitem{DR2} G. Da Prato and M. R\"ockner, \emph{Invariant measures for a
stochastic porous medium equation,} preprint SNS, 2003; to appear
in Proceedings of Conference in Honour of K. It\^o, Kyoto, 2002.

\bibitem{DRRW} G. Da Prato, M. R\"ockner, Rozovskii and F.-Y.
Wang, \emph{Strong solutions to stochastic generalized porous
media equations: existence, uniqueness and ergodicity,} to appear
in Comm. Part. Diff. Equat.

\bibitem{DaZa} G. Da Prato and J. Zabczyk, \emph{ Stochastic Equations in
Infinite Dimensions,} Encyclopedia of Mathematics and its
Applications, Cambridge University Press. 1992.

\bibitem{DZ} A. Dembo and O. Zeitouni, \emph{ Large deviations
Techniques and Applications. Second Edition,} Springer, New York.
1998.

\bibitem{FW} M.I. Freidlin and A.D. Wentzell, \emph{Random perturbations of dynamical systems,} Translated from the Russian by Joseph Szu"cs. Grundlehren der Mathematischen Wissenschaften [Fundamental Principles of Mathematical Sciences], 260. Springer-Verlag, New York, 1984.

\bibitem{Hino} M. Hino, J.A. Rami'rez, \emph{Small-time Gaussian behavior of symmetric diffusion semigroups,}
  Ann. Probab.  31  (2003),  no. 3, 1254--1295.

\bibitem{KR1} N.V. Krylov and M. R\"ockner, \emph{Strong solutions of stochastic
equations with singular time dependent drift,} to appear in
Probab. Theory Relat. Fields.

\bibitem{KR} N.V. Krylov and B.L. Rozovskii, \emph{Stochastic evolution equations,}
Translated from Itogi Naukii Tekhniki, Seriya Sovremennye Problemy
Matematiki 14(1979), 71--146, Plenum Publishing Corp. 1981.

\bibitem{ma-roe92}
  Z.M. Ma and M. R\"ockner,
  \emph{Introduction to the theory of (non-symmetric) Dirichlet forms},
  Springer, 1992.

\bibitem{RRW} J. Ren, M. R\"ockner and F.-Y. Wang,
\emph{Stochastic generalized porous media and fast diffusion
equations,} Preprint 2005.

\bibitem{Var} S.R.S. Varadhan, \emph{Diffusion processes in a small time interval,}
  Comm. Pure Appl. Math.  20  (1967), 659--685.

\bibitem{Wu04} L. Wu,
\emph{On large deviations for moving average processes, } In
Probability, Finance and Insurance, pp.15-49, the proceeding of a
Workshop at the University of Hong-Kong (15-17 July 2002), Eds:
T.L. Lai, H.L. Yang and S.P. Yung. World Scientific 2004,
Singapour.

\bibitem{Z} E. Zeidler, \emph{Nonlinear Functional Analysis and
its Applications, II/B, Nonlinear Monotone Operators,}
Springer-Verlag, New York: 1990.

\end{thebibliography}

\end{document}